\documentclass[11pt,english]{article}
\usepackage[T1]{fontenc}
\usepackage[latin9]{inputenc}
\usepackage{geometry}
\usepackage{refstyle}
\usepackage{amsthm}
\usepackage{amsmath}
\usepackage{amssymb}
\usepackage{setspace}
\usepackage{esint}
\makeatletter

\AtBeginDocument{\providecommand\secref[1]{\ref{sec:#1}}}
\AtBeginDocument{\providecommand\eqref[1]{\ref{eq:#1}}}
\AtBeginDocument{\providecommand\lemref[1]{\ref{lem:#1}}}
\AtBeginDocument{\providecommand\propref[1]{\ref{prop:#1}}}
\RS@ifundefined{subref}
  {\def\RSsubtxt{section~}\newref{sub}{name = \RSsubtxt}}
  {}
\RS@ifundefined{thmref}
  {\def\RSthmtxt{theorem~}\newref{thm}{name = \RSthmtxt}}
  {}
\RS@ifundefined{lemref}
  {\def\RSlemtxt{lemma~}\newref{lem}{name = \RSlemtxt}}
  {}

\theoremstyle{plain}
\newtheorem{thm}{\protect\theoremname}
  \theoremstyle{plain}
  \newtheorem{lem}[thm]{\protect\lemmaname}
  \theoremstyle{plain}
  \newtheorem{prop}[thm]{\protect\propositionname}

\newref{eq}{name = {}, refcmd = \textup{(\ref{#1})}}
\newref{prop}{name = Proposition~}
\newref{lem}{name = Lemma~}
\newref{sec}{name = Sect.~}

\DeclareMathSymbol{\Gamma}{\mathalpha}{operators}{0}
\DeclareMathSymbol{\Delta}{\mathalpha}{operators}{1}
\DeclareMathSymbol{\Theta}{\mathalpha}{operators}{2}
\DeclareMathSymbol{\Lambda}{\mathalpha}{operators}{3}
\DeclareMathSymbol{\Xi}{\mathalpha}{operators}{4}
\DeclareMathSymbol{\Pi}{\mathalpha}{operators}{5}
\DeclareMathSymbol{\Sigma}{\mathalpha}{operators}{6}
\DeclareMathSymbol{\Upsilon}{\mathalpha}{operators}{7}
\DeclareMathSymbol{\Phi}{\mathalpha}{operators}{8}
\DeclareMathSymbol{\Psi}{\mathalpha}{operators}{9}
\DeclareMathSymbol{\Omega}{\mathalpha}{operators}{10}

\makeatother

\usepackage{babel}
  \providecommand{\lemmaname}{Lemma}
  \providecommand{\propositionname}{Proposition}
\providecommand{\theoremname}{Theorem}

\begin{document}

\title{The $\Delta_{2}$-condition and $\varphi$-families\\
of probability distributions}

\author{Rui F.\ Vigelis%
\thanks{Computer Engineering, Campus Sobral, Federal University of Ceará,
Sobral-CE, Brazil. E-mail: rfvigelis@ufc.br.%
},\and Charles C.\ Cavalcante%
\thanks{Wireless Telecommunication Research Group, Department of Teleinformatics
Engineering, Federal University of Ceará, Fortaleza-CE, Brazil. E-mail:
charles@ufc.br.%
}}
\maketitle
\begin{abstract}
In this paper, we provide some results related to the $\Delta_{2}$-condition
of Musielak--Orlicz functions and $\varphi$-families of probability
distributions, which are modeled on Musielak--Orlicz spaces. We show
that if two $\varphi$-families are modeled on Musielak--Orlicz spaces
generated by Musielak--Orlicz functions satisfying the $\Delta_{2}$-condition,
then these $\varphi$-families are equal as sets. We also investigate
the behavior of the normalizing function near the boundary of the
set on which a $\varphi$-family is defined.
\end{abstract}

\section{Introduction}

In \cite{Vigelis:2011}, $\varphi$-families of probability distributions
are introduced as a generalization of exponential families of probability
distributions \cite{Pistone:1995,Pistone:1999}. The main idea leading
to this generalization is the replacement of the exponential function
with a $\varphi$-function (a definition is given below). These families
(of probability distributions) are subsets of the collection $\mathcal{P}_{\mu}$
of all $\mu$-a.e.\ strictly positive probability densities. What
the papers \cite{Pistone:1995,Pistone:1999,Vigelis:2011} provide
is a framework endowing $\mathcal{P}_{\mu}$ with a structure of $C^{\infty}$-Banach
manifold \cite{Lang:1995}, where a family constitutes a connected
component of $\mathcal{P}_{\mu}$. These families are modeled on Musielak--Orlicz
spaces (exponential families are modeled on exponential Orlicz spaces)
\cite{Musielak:1983,Krasnoselskii:1961,Rao:1991}. In many properties
of these spaces, the $\Delta_{2}$-condition of Musielak--Orlicz functions
plays a central role. For example, a Musielak--Orlicz space $L^{\Phi}$
is equal to the Musielak--Orlicz class $\tilde{L}^{\Phi}$ if and
only if the Musielak--Orlicz function $\Phi$ satisfies the $\Delta_{2}$-condition.
In this paper we investigate the $\Delta_{2}$-condition in the context
of $\varphi$-families. In \secref{Delta2_phi-families}, we show
that if two $\varphi$-families are modeled on Musielak--Orlicz spaces
generated by Musielak--Orlicz functions satisfying the $\Delta_{2}$-condition,
then these $\varphi$-families are equal as sets. In \secref{psi_behavior},
we investigate the behavior of the normalizing function near the boundary
of the set on which a $\varphi$-family is defined. In the rest of
this section, $\varphi$-families are exposed.

A $\varphi$-family is the image of a mapping whose domain is a subset
of a Musielak--Orlicz space. In what follows, this statement will
be made more precise. Musielak--Orlicz spaces are just briefly introduced
here. These spaces are thoroughly exposed in \cite{Musielak:1983,Krasnoselskii:1961,Rao:1991}. 

Let $(T,\Sigma,\mu)$ be a $\sigma$-finite, non-atomic measure space.
A function $\Phi\colon T\times[0,\infty)\rightarrow[0,\infty]$ is
said to be a \textit{Musielak--Orlicz function} if
\begin{itemize}
\item [(i)] $\Phi(t,\cdot)$ is convex and lower semi-continuous for $\mu$-a.e.\ $t\in T$,
\item [(ii)] $\Phi(t,0)=\lim_{u\downarrow0}\Phi(t,u)=0$ and $\lim_{u\rightarrow\infty}\Phi(t,u)=\infty$
for $\mu$-a.e.\ $t\in T$,
\item [(iii)] $\Phi(\cdot,u)$ is measurable for each $u\geq0$.
\end{itemize}
We notice that $\Phi(t,\cdot)$, by (i)--(ii), is not equal to $0$
or $\infty$ on the interval $(0,\infty)$. A Musielak--Orlicz function
$\Phi$ is said to be an \textit{Orlicz function} if the functions
$\Phi(t,\cdot)$ are the same for $\mu$-a.e.\ $t\in T$.

Let $L^{0}$ denote the linear space of all real-valued, measurable
functions on $T$, with equality $\mu$-a.e. Given any Musielak--Orlicz
function $\Phi$, we denote the functional $I_{\Phi}(u)=\int_{T}\Phi(t,|u(t)|)d\mu$,
for any $u\in L^{0}$. The \textit{Musielak--Orlicz space}, \textit{Musielak--Orlicz
class}, and \textit{Morse--Transue space} generated by a Musielak--Orlicz
function $\Phi$ are defined by 
\begin{align*}
L^{\Phi} & =\{u\in L^{0}:I_{\Phi}(\lambda u)<\infty\text{ for some }\lambda>0\},\\
\tilde{L}^{\Phi} & =\{u\in L^{0}:I_{\Phi}(u)<\infty\},\\
\intertext{and}E^{\Phi} & =\{u\in L^{0}:I_{\Phi}(\lambda u)<\infty\text{ for all }\lambda>0\},
\end{align*}
respectively. The Musielak--Orlicz space $L^{\Phi}$ is a Banach space
when it is equipped with the \textit{Luxemburg norm} 
\[
\Vert u\Vert_{\Phi}=\inf\Bigl\{\lambda>0:I_{\Phi}\Bigl(\frac{u}{\lambda}\Bigr)\leq1\Bigr\},
\]
or the \textit{Orlicz norm}
\[
\Vert u\Vert_{\Phi,0}=\sup\biggl\{\biggl|\int_{T}uvd\mu\biggr|:v\in\tilde{L}^{\Phi^{*}}\text{ and }I_{\Phi^{*}}(v)\leq1\biggr\},
\]
where $\Phi^{*}(t,v)=\sup\nolimits _{u\geq0}(uv-\Phi(t,u))$ is the
\textit{Fenchel conjugate} of $\Phi(t,\cdot)$. These norms are equivalent
and the inequalities $\Vert u\Vert_{\Phi}\leq\Vert u\Vert_{\Phi,0}\leq2\Vert u\Vert_{\Phi}$
hold for all $u\in L^{\Phi}$.

Whereas exponential families are based on the exponential function,
$\varphi$-families are based on $\varphi$-functions. A function
$\varphi\colon T\times\mathbb{R}\rightarrow(0,\infty)$ is said to
be a \textit{$\varphi$-function} if the following conditions are
satisfied:
\begin{itemize}
\item [(a1)] $\varphi(t,\cdot)$ is convex for $\mu$-a.e.\ $t\in T$,
\item [(a2)] $\lim_{u\rightarrow-\infty}\varphi(t,u)=0$ and $\lim_{u\rightarrow\infty}\varphi(t,u)=\infty$
for $\mu$-a.e.\ $t\in T$,
\item [(a3)] $\varphi(\cdot,u)$ is measurable for each $u\in\mathbb{R}$.
\end{itemize}
In addition, we assume a positive, measurable function $u_{0}\colon T\rightarrow(0,\infty)$
can be found such that, for every measurable function $c\colon T\rightarrow\mathbb{R}$
for which $\varphi(t,c(t))$ is in $\mathcal{P}_{\mu}$, we have that
\begin{itemize}
\item [(a4)] $\varphi(t,c(t)+\lambda u_{0}(t))$ is $\mu$-integrable for
all $\lambda>0$.
\end{itemize}
The exponential function is an example of $\varphi$-function, since
$\varphi(t,u)=\exp(u)$ satisfies conditions (a1)--(a3) and (a4) with
$u_{0}=\mathbf{1}_{T}$, where $\boldsymbol{1}_{A}$ is the indicator
function of a subset $A\subseteq T$. Another example of $\varphi$-function
is the Kaniadakis' $\kappa$-exponential (see \cite{Kaniadakis:2002}
and \cite[Example~1]{Vigelis:2011}). Let $\varphi_{+}'(t,\cdot)$
denote the right derivative of $\varphi(t,\cdot)$. In what follows,
$\boldsymbol{\varphi}$ and $\boldsymbol{\varphi}_{+}'$ denote the
function operators $\boldsymbol{\varphi}(u)(t):=\varphi(t,u(t))$
and $\boldsymbol{\varphi}_{+}'(u)(t):=\varphi_{+}'(t,u(t))$, respectively,
for any real-valued function $u\colon T\rightarrow\mathbb{R}$. 

A $\varphi$-family is defined to be a subset of the collection
\[
\mathcal{P}_{\mu}=\{p\in L^{0}:p>0\text{ and }\mathbb{E}[p]=1\},
\]
where $\mathbb{E}[\cdot]=\int_{T}(\cdot)d\mu$ denotes integration
with respect to $\mu$. For each probability density $p\in\mathcal{P}_{\mu}$,
we associate a $\varphi$-family $\mathcal{F}_{c}^{\varphi}\subset\mathcal{P}_{\mu}$
centered at $p$, where $c\colon T\rightarrow\mathbb{R}$ is a measurable
function such that $p=\boldsymbol{\varphi}(c)$. The Musielak--Orlicz
space $L^{\Phi_{c}}$ on which the $\varphi$-family $\mathcal{F}_{c}^{\varphi}$
is modeled is given in terms of the Musielak--Orlicz function 
\begin{equation}
\Phi_{c}(t,u)=\varphi(t,c(t)+u)-\varphi(t,c(t)).\label{eq:MO_phi_c}
\end{equation}
We will use the notation $L_{c}^{\varphi}$, $\tilde{L}_{c}^{\varphi}$
and $E_{c}^{\varphi}$ in the place of $L^{\Phi_{c}}$, $\tilde{L}^{\Phi_{c}}$
and $E^{\Phi_{c}}$, respectively, to indicate that $\Phi_{c}$ is
given by \eqref{MO_phi_c}. Because $\boldsymbol{\varphi}(c)$ is
$\mu$-integrable, the Musielak--Orlicz space $L_{c}^{\varphi}$ corresponds
to the set of all functions $u\in L^{0}$ for which there exists $\varepsilon>0$
such that $\boldsymbol{\varphi}(c+\lambda u)$ is $\mu$-integrable
for all $\lambda\in(-\varepsilon,\varepsilon)$.

The elements of the $\varphi$-family $\mathcal{F}_{c}^{\varphi}\subset\mathcal{P}_{\mu}$
centered at $p=\boldsymbol{\varphi}(c)\in\mathcal{P}_{\mu}$ are given
by the one-to-one mapping 
\begin{equation}
\boldsymbol{\varphi}_{c}(u):=\boldsymbol{\varphi}(c+u-\psi(u)u_{0}),\qquad\text{for each }u\in\mathcal{B}_{c}^{\varphi},\label{eq:phi_c-1-1}
\end{equation}
where the set $\mathcal{B}_{c}^{\varphi}\subseteq L_{c}^{\varphi}$
is defined as the intersection of the convex set
\[
\mathcal{K}_{c}^{\varphi}=\{u\in L_{c}^{\varphi}:\mathbb{E}[\boldsymbol{\varphi}(c+\lambda u)]<\infty\text{ for some }\lambda>1\}
\]
with the closed subspace
\[
B_{c}^{\varphi}=\{u\in L_{c}^{\varphi}:\mathbb{E}[u\boldsymbol{\varphi}_{+}'(c)]=0\},
\]
and the \textit{normalizing function} $\psi\colon\mathcal{B}_{c}^{\varphi}\rightarrow[0,\infty)$
is introduced so that expression \eqref{phi_c-1-1} defines a probability
distribution in $\mathcal{P}_{\mu}$. By \cite[Lemma~2]{Vigelis:2011},
the set $\mathcal{K}_{c}^{\varphi}$ is open in $L_{c}^{\varphi}$,
and hence $\mathcal{B}_{c}^{\varphi}$ is open in $B_{c}^{\varphi}$. 

Its is clear that the collection $\{\mathcal{F}_{c}^{\varphi}:\boldsymbol{\varphi}(c)\in\mathcal{P}_{\mu}\}$
covers the whole family $\mathcal{P}_{\mu}$. Moreover, $\varphi$-families
are maximal in the sense that if two $\varphi$-families have a non-empty
intersection, then they coincide as sets. Let $\mathcal{F}_{c_{1}}^{\varphi}$
and $\mathcal{F}_{c_{2}}^{\varphi}$ be two $\varphi$-families centered
at $\boldsymbol{\varphi}(c_{1})\in\mathcal{P}_{\mu}$ and $\boldsymbol{\varphi}(c_{2})\in\mathcal{P}_{\mu}$,
for some measurable functions $c_{1},c_{2}\colon T\rightarrow\mathbb{R}$.
If the $\varphi$-families $\mathcal{F}_{c_{1}}^{\varphi}$ and $\mathcal{F}_{c_{2}}^{\varphi}$
have non-empty intersection, then $\mathcal{F}_{c_{1}}^{\varphi}=\mathcal{F}_{c_{2}}^{\varphi}$
and the spaces $L_{c_{1}}^{\varphi}$ and $L_{c_{2}}^{\varphi}$ are
equal as sets, and have equivalent norms. Because the transition map
$\boldsymbol{\varphi}_{c_{2}}^{-1}\circ\boldsymbol{\varphi}_{c_{1}}\colon\mathcal{B}_{c_{1}}^{\varphi}\rightarrow\mathcal{B}_{c_{2}}^{\varphi}$
is an affine transformation, the collection of charts $\{(\mathcal{B}_{c}^{\varphi},\boldsymbol{\varphi}_{c})\}_{\boldsymbol{\varphi}(c)\in\mathcal{P}_{\mu}}$
is an atlas of class $C^{\infty}$, endowing $\mathcal{P}_{\mu}$
with a structure of $C^{\infty}$-Banach manifold. A verification
of these claims is found in \cite{Vigelis:2011}.

\section{The $\Delta_{2}$-condition and $\varphi$-families\label{sec:Delta2_phi-families}}

A Musielak--Orlicz function $\Phi$ is said to satisfy the \textit{$\Delta_{2}$-condition},
or to belong to the \textit{$\Delta_{2}$-class} (denoted by $\Phi\in\Delta_{2}$),
if a constant $K>0$ and a non-negative function $f\in\tilde{L}^{\Phi}$
can be found such that
\begin{equation}
\Phi(t,2u)\leq K\Phi(t,u),\qquad\text{for all }u\geq f(t),\quad\text{and }\mu\text{-a.e. }t\in T.\label{eq:D2_0}
\end{equation}
It is easy to see that, if a Musielak--Orlicz function $\Phi$ satisfies
the $\Delta_{2}$-condition, then $I_{\Phi}(u)<\infty$ for every
$u\in L^{\Phi}$. In this case, $L^{\Phi}$, $\tilde{L}^{\Phi}$ and
$E^{\Phi}$ are equal as sets. On the other hand, if the Musielak--Orlicz
function $\Phi$ does not satisfy the $\Delta_{2}$-condition, then
$E^{\Phi}$ is a proper subspace of $L^{\Phi}$. In addition, we can
state:

\begin{lem}
\label{lem:u_stars} Let $\Phi$ be a Musielak--Orlicz function not
satisfying the $\Delta_{2}$-condition and such that $\Phi(t,b_{\Phi}(t))=\infty$
for $\mu$-a.e.\ $t\in T$, where $b_{\Phi}(t)=\sup\{u\geq0:\Phi(t,u)<\infty\}$.
Then we can find functions $u_{*}$ and $u^{*}$ in $L^{\Phi}$ such
that 
\begin{equation}
\left\{ \begin{array}{ll}
I_{\Phi}(\lambda u_{*})<\infty, & \quad\text{for }0\leq\lambda\leq1,\\
I_{\Phi}(\lambda u_{*})=\infty, & \quad\text{for }1<\lambda,
\end{array}\right.\label{eq:substar}
\end{equation}
and
\begin{equation}
\left\{ \begin{array}{ll}
I_{\Phi}(\lambda u^{*})<\infty, & \quad\text{for }0\leq\lambda<1,\\
I_{\Phi}(\lambda u^{*})=\infty, & \quad\text{for }1\leq\lambda.
\end{array}\right.\label{eq:upstar}
\end{equation}

\end{lem}

This lemma is a well established result for Orlicz functions (see
\cite[Sect.~8.4]{Krasnoselskii:1961}). A proof of \lemref{u_stars}
is given in \cite{Vigelis:2013}. The next result shows that we can
always find a $\varphi$-family modeled on a Musielak--Orlicz space
generated by a Musielak--Orlicz function not satisfying the $\Delta_{2}$-condition.

\begin{prop}
Given any $\varphi$-function $\varphi$, we can find a measurable
function $c\colon T\rightarrow\mathbb{R}$ with $\mathbb{E}[\boldsymbol{\varphi}(c)]=1$
such that the Musielak--Orlicz function $\Phi_{c}(t,u)=\varphi(t,c(t)+u)-\varphi(t,c(t))$
does not satisfy the $\Delta_{2}$-condition.\end{prop}
\begin{proof}
Let $A$ and $B$ be two disjoint, measurable sets satisfying $0<\mu(A)<\infty$
and $0<\mu(B)<\infty$. Fixed any measurable function $\widetilde{c}$
such that $\mathbb{E}[\boldsymbol{\varphi}(\widetilde{c})]=1$, we
take any non-integrable function $f$ supported on $A$ such that
$\boldsymbol{\varphi}(\widetilde{c})\mathbf{1}_{A}\leq f\mathbf{1}_{A}<\infty$.
Let $u\colon T\rightarrow[0,\infty)$ be a measurable function supported
on $A$ such that $\boldsymbol{\varphi}(\widetilde{c}+u)\mathbf{1}_{A}=f\mathbf{1}_{A}$.
If $\beta>0$ is such that $\mathbb{E}[\boldsymbol{\varphi}(\widetilde{c}-u)\mathbf{1}_{A}]+\beta\mu(B)+\mathbb{E}[\boldsymbol{\varphi}(\widetilde{c})\mathbf{1}_{T\setminus(A\cup B)}]=1$,
then we define
\[
c=(\widetilde{c}-u)\mathbf{1}_{A}+\overline{c}\mathbf{1}_{B}+\widetilde{c}\mathbf{1}_{T\setminus(A\cup B)},
\]
where $\overline{c}\colon T\rightarrow\mathbb{R}$ is a measurable
function supported on $B$ such that $\varphi(t,\overline{c}(t))=\beta$,
for $\mu$-a.e.\ $t\in B$. Because the function $u$ is supported
on $A$, we can write
\[
\mathbb{E}[\boldsymbol{\varphi}(c+u)]=\mathbb{E}[\boldsymbol{\varphi}(\widetilde{c})\mathbf{1}_{A}]+\mathbb{E}[\boldsymbol{\varphi}(\overline{c})\mathbf{1}_{B}]+\mathbb{E}[\boldsymbol{\varphi}(\widetilde{c})\mathbf{1}_{T\setminus(A\cup B)}]<\infty.
\]
On the other hand, since $f$ is non-integrable, we have
\[
\mathbb{E}[\boldsymbol{\varphi}(c+2u)]>\mathbb{E}[\boldsymbol{\varphi}(\widetilde{c}+u)\mathbf{1}_{A}]=\mathbb{E}[f]=\infty.
\]
Therefore, the Musielak--Orlicz function $\Phi_{c}$ does not satisfy
the $\Delta_{2}$-condition.
\end{proof}

The main result of this section is a consequence of the following
proposition:

\begin{prop}
\label{prop:Lvarphib_subset_Lvarphic} Let $b:T\rightarrow\mathbb{R}$
be a measurable function such that $\mathbb{E}[\boldsymbol{\varphi}(b)]=1$.
Then $L_{b}^{\varphi}\subseteq L_{c}^{\varphi}$ for every measurable
function $c:T\rightarrow\mathbb{R}$ such that $\mathbb{E}[\boldsymbol{\varphi}(c)]=1$
if, and only if, the Musielak--Orlicz function $\Phi_{b}(t,u)=\varphi(t,b(t)+u)-\varphi(t,b(t))$
satisfies the $\Delta_{2}$-condition.\end{prop}
\begin{proof}
Assume that $\Phi_{b}$ satisfies the $\Delta_{2}$-condition. Let
$c:T\rightarrow\mathbb{R}$ be any measurable function such that $\mathbb{E}[\boldsymbol{\varphi}(c)]=1$.
Denoting $A=\{t\in T:c(t)\geq b(t)\}$, it is clear that the function
$(c-b)\mathbf{1}_{A}$ is in $L_{b}^{\varphi}$. Hence, for any function
$u\in L_{b}^{\varphi}$, we can write
\[
\mathbb{E}[\boldsymbol{\varphi}(c+|u|)]=\mathbb{E}[\boldsymbol{\varphi}(b+(c-b)+|u|)]\leq\mathbb{E}[\boldsymbol{\varphi}(b+(c-b)\mathbf{1}_{A}+|u|)]<\infty,
\]
since $(c-b)\mathbf{1}_{A}+|u|$ is in $L_{b}^{\varphi}$, and the
sets $L_{b}^{\varphi}$ and $\tilde{L}_{b}^{\varphi}$ are equal.
Thus, $L_{b}^{\varphi}\subseteq L_{c}^{\varphi}$.

Now we suppose that $\Phi_{b}$ does not satisfy the $\Delta_{2}$-condition.
From \lemref{u_stars}, there exists a non-negative function $u\in\tilde{L}^{\Phi_{b}}$
such that $I_{\Phi_{b}}(\lambda u)=\infty$ for all $\lambda>1$.
Using the function $u$, we will provide a measurable function $c:T\rightarrow\mathbb{R}$
with $\mathbb{E}[\boldsymbol{\varphi}(c)]=1$ for which $L_{b}^{\varphi}$
is not contained in $L_{c}^{\varphi}$. By \cite{Kaminska:1985} or
\cite[Lemma~2]{Kolwicz:1995}, we can find a sequence of non-decreasing,
measurable sets $\{T_{n}\}$, satisfying $\mu(T_{n})<\infty$ and
$\mu(T\setminus\bigcup_{n=1}^{\infty}T_{n})=0$, such that
\begin{equation}
\operatorname*{ess\, sup}_{t\in T_{n}}\Phi_{b}(t,u)<\infty,\qquad\text{for all }u>0,\text{ and each }n\geq1.\label{eq:sup_T_n-1-1}
\end{equation}
Thus, for a sufficiently large $n_{0}\geq1$, the set $A=\{t\in T_{n_{0}}:u(t)\leq n_{0}\}$
satisfies $\mathbb{E}[\boldsymbol{\varphi}(b+u)\mathbf{1}_{T\setminus A}]<1$.
Observing that
\[
I_{\Phi_{b}}(\lambda u\mathbf{1}_{A})\leq\Bigl[\operatorname*{ess\, sup}_{t\in T_{n_{0}}}\Phi_{b}(t,\lambda n_{0})\Bigr]\mu(T_{n_{0}})<\infty,\qquad\text{for each }\lambda>0,
\]
we can infer that
\begin{equation}
I_{\Phi_{b}}(\lambda u\mathbf{1}_{T\setminus A})=I_{\Phi_{b}}(\lambda u)-I_{\Phi_{b}}(\lambda u\mathbf{1}_{A})=\infty,\qquad\text{for all }\lambda>1.\label{eq:IPhi_infty}
\end{equation}
Let $\alpha>0$ be such that $\alpha\mu(A)+\mathbb{E}[\boldsymbol{\varphi}(b+u)\mathbf{1}_{T\setminus A}]=1$.
Then we define
\[
c=\overline{c}\mathbf{1}_{A}+(b+u)\mathbf{1}_{T\setminus A},
\]
where $\overline{c}\colon T\rightarrow\mathbb{R}$ is a measurable
function supported on $A$ such that $\varphi(t,\overline{c}(t))=\alpha$,
for $\mu$-a.e.\ $t\in A$. It is clear that $\mathbb{E}[\boldsymbol{\varphi}(c)]=1$.
According to \cite[Proposition~4]{Vigelis:2011}, if $c_{1},c_{2}\colon T\rightarrow\mathbb{R}$
are measurable functions such that $\mathbb{E}[\boldsymbol{\varphi}(c_{1})]=1$
and $\mathbb{E}[\boldsymbol{\varphi}(c_{2})]=1$, then $(c_{1}-c_{2})\in L_{c_{2}}^{\varphi}$
is a necessary and sufficient condition for $L_{c_{1}}^{\varphi}\subseteq L_{c_{2}}^{\varphi}$.
Thus, to show that $L_{b}^{\varphi}$ is not contained in $L_{c}^{\varphi}$,
we have to verify that $(b-c)\notin L_{c}^{\varphi}$. Denoting $F=\{t\in T:c(t)\geq b(t)\}$,
for any $\lambda>0$, we can write
\begin{align}
\mathbb{E}[\boldsymbol{\varphi}(c+\lambda|b-c|)] & \geq\mathbb{E}[\boldsymbol{\varphi}(c+\lambda(c-b))\mathbf{1}_{F}]\nonumber \\
 & =\mathbb{E}[\boldsymbol{\varphi}(b+(1+\lambda)(c-b))\mathbf{1}_{F}]\nonumber \\
 & \geq\mathbb{E}[\boldsymbol{\varphi}(b+(1+\lambda)u)\mathbf{1}_{T\setminus A}]\label{eq:Ephi_ineq}\\
 & =\infty,\label{eq:Ephi_infty}
\end{align}
where in \eqref{Ephi_ineq} we used that $T\setminus A\subseteq F$
and $(c-b)\mathbf{1}_{T\setminus A}=u\mathbf{1}_{T\setminus A}$,
and \eqref{Ephi_infty} follows from \eqref{IPhi_infty}. We conclude
that $(b-c)\notin L_{c}^{\varphi}$, and hence $L_{b}^{\varphi}$
is not contained in $L_{c}^{\varphi}$. Therefore, if $L_{b}^{\varphi}\subseteq L_{c}^{\varphi}$
for any measurable function $c:T\rightarrow\mathbb{R}$ such that
$\mathbb{E}[\boldsymbol{\varphi}(c)]=1$, then the Musielak--Orlicz
function $\Phi_{b}$ satisfies the $\Delta_{2}$-condition.
\end{proof}

Now we can state the main result of this section:

\begin{prop}
Let $b,c\colon T\rightarrow\mathbb{R}$ be measurable functions such
that $\mathbb{E}[\boldsymbol{\varphi}(b)]=1$ and $\mathbb{E}[\boldsymbol{\varphi}(c)]=1$.
If the Musielak--Orlicz functions $\Phi_{b}(t,u)=\varphi(t,b(t)+u)-\varphi(t,b(t))$
and $\Phi_{c}(t,u)=\varphi(t,c(t)+u)-\varphi(t,c(t))$ satisfy the
$\Delta_{2}$-condition, then $L_{b}^{\varphi}$ and $L_{c}^{\varphi}$
are equal as sets. Moreover, $\mathcal{F}_{b}^{\varphi}=\mathcal{F}_{c}^{\varphi}$.\end{prop}
\begin{proof}
The conclusion that $L_{b}^{\varphi}$ and $L_{c}^{\varphi}$ are
equal as sets follows from \propref{Lvarphib_subset_Lvarphic}. By
\cite[Proposition~4]{Vigelis:2011}, it is clear that $(c-b)\in\mathcal{K}_{b}^{\varphi}$.
Let $\alpha\geq0$ be such that $u=(c-b)+\alpha u_{0}$ belongs to
$\mathcal{B}_{b}^{\varphi}$. If $\psi_{1}$ is the normalizing function
associated with $\mathcal{F}_{b}^{\varphi}$, then $\psi_{1}(u)=\alpha$
and $\boldsymbol{\varphi}_{b}(u)=\boldsymbol{\varphi}(b+u-\psi_{1}(u)u_{0})=\boldsymbol{\varphi}(c)$.
Thus the $\varphi$-families $\mathcal{F}_{b}^{\varphi}$ and $\mathcal{F}_{c}^{\varphi}$
have a non-empty intersection, and hence $\mathcal{F}_{b}^{\varphi}=\mathcal{F}_{c}^{\varphi}$.
\end{proof}

\section{The behavior of $\psi$ near the boundary of $\mathcal{B}_{c}^{\varphi}$\label{sec:psi_behavior}}

In this section, we investigate the behavior of the normalizing function
$\psi$ near the boundary of $\mathcal{B}_{c}^{\varphi}$ (with respect
to the topology of $B_{c}^{\varphi}$). More specifically, given any
function $u$ in the boundary of $\mathcal{B}_{c}^{\varphi}$, which
we denote by $\partial\mathcal{B}_{c}^{\varphi}$, we want to know
whether $\psi(\lambda u)$ converges to a finite value or not as $\lambda\uparrow1$.
For this purpose, we establish under what conditions the set $\mathcal{B}_{c}^{\varphi}$
has a non-empty boundary. This result is related to the $\Delta_{2}$-condition.
By definition, a function $u\in L^{0}$ is in $\mathcal{K}_{c}^{\varphi}$
if there exists $\varepsilon>0$ such that $\mathbb{E}[\boldsymbol{\varphi}(c+\lambda u)]<\infty$
for all $\lambda\in(-\varepsilon,1+\varepsilon)$. Because the set
$\mathcal{B}_{c}^{\varphi}=\mathcal{K}_{c}^{\varphi}\cap B_{c}^{\varphi}$
is open in $B_{c}^{\varphi}$, we conclude that a function $u\in B_{c}^{\varphi}$
belongs to the boundary of $\mathcal{B}_{c}^{\varphi}$ if and only
if $\mathbb{E}[\boldsymbol{\varphi}(c+\lambda u)]<\infty$ for all
$\lambda\in(0,1)$, and $\mathbb{E}[\boldsymbol{\varphi}(c+\lambda u)]=\infty$
for each $\lambda>1$. If the Musielak--Orlicz function $\Phi_{c}=\varphi(t,c(t)+u)-\varphi(t,c(t))$
satisfies the $\Delta_{2}$-condition, then $\mathbb{E}[\boldsymbol{\varphi}(c+u)]<\infty$
for all $u\in L_{c}^{\varphi}$. In this case, the set $\mathcal{B}_{c}^{\varphi}$
coincides with the closed subspace $B_{c}^{\varphi}$, and the boundary
of $\mathcal{B}_{c}^{\varphi}$ is empty. On the other hand, if $\Phi_{c}$
does not satisfies the $\Delta_{2}$-condition, then the boundary
of $\mathcal{B}_{c}^{\varphi}$ is non-empty. Moreover, not all functions
$u$ in the boundary of $\mathcal{B}_{c}^{\varphi}$ satisfy $\mathbb{E}[\boldsymbol{\varphi}(c+u)]<\infty$
(or $\mathbb{E}[\boldsymbol{\varphi}(c+u)]=\infty$). In other words,
we can always find functions $w_{*}$ and $w^{*}$ in $\partial\mathcal{B}_{c}^{\varphi}$
for which $\mathbb{E}[\boldsymbol{\varphi}(c+w_{*})]<\infty$ and
$\mathbb{E}[\boldsymbol{\varphi}(c+w^{*})]=\infty$. This result,
which is a consequence of \lemref{u_stars}, is provided by the following
proposition:

\begin{prop}
\label{prop:w_stars} The boundary of $\mathcal{B}_{c}^{\varphi}$
is non-empty if and only if the Musielak--Orlicz function $\Phi_{c}=\varphi(t,c(t)+u)-\varphi(t,c(t))$
does not satisfy the $\Delta_{2}$-condition. Moreover, in any of
these cases, there exist functions $w_{*}$ and $w^{*}$ in $\partial\mathcal{B}_{c}^{\varphi}$
such that $\mathbb{E}[\boldsymbol{\varphi}(c+w_{*})]<\infty$ and
\textup{$\mathbb{E}[\boldsymbol{\varphi}(c+w^{*})]=\infty$.}\end{prop}
\begin{proof}
Given non-negative functions $u_{*}$ and $u^{*}$ in $L_{c}^{\varphi}$
satisfying \eqref{substar} and \eqref{upstar} in \lemref{u_stars},
we consider the functions
\[
w_{*}=u_{*}-\frac{\mathbb{E}[u_{*}\boldsymbol{\varphi}_{+}'(c)]}{\mathbb{E}[u_{0}\boldsymbol{\varphi}_{+}'(c)]}u_{0},\qquad\text{and}\qquad w^{*}=u^{*}-\frac{\mathbb{E}[u^{*}\boldsymbol{\varphi}_{+}'(c)]}{\mathbb{E}[u_{0}\boldsymbol{\varphi}_{+}'(c)]}u_{0},
\]
which are in $B_{c}^{\varphi}$. Next we show that $w_{*}$ is in
$\partial\mathcal{B}_{c}^{\varphi}$ and satisfies $\mathbb{E}[\boldsymbol{\varphi}(c+w_{*})]<\infty$.
For any $0\leq\lambda\leq1$, its clear that 
\[
\mathbb{E}[\boldsymbol{\varphi}(c+\lambda w_{*})]\leq\mathbb{E}[\boldsymbol{\varphi}(c+\lambda u_{*})]<\infty.
\]
Now suppose that $\mathbb{E}[\boldsymbol{\varphi}(c+\lambda_{0}w_{*})]<\infty$
for some $\lambda_{0}>1$. In view of $1\leq\mathbb{E}[\boldsymbol{\varphi}(c+\lambda_{0}w_{*})]<\infty$,
we can find $\alpha_{0}\geq0$ such that $\mathbb{E}[\boldsymbol{\varphi}(c+\lambda_{0}w_{*}-\alpha_{0}u_{0})]=1$.
By the definition of $u_{0}$, fixed any measurable function $\widetilde{c}$
such that $\mathbb{E}[\boldsymbol{\varphi}(\widetilde{c})]=1$, we
have that $\mathbb{E}[\boldsymbol{\varphi}(\widetilde{c}+\alpha u_{0})]<\infty$
for all $\alpha\in\mathbb{R}$. Hence, considering $\widetilde{c}=c+\lambda_{0}w_{*}-\alpha_{0}u_{0}$
and 
\[
\alpha=\lambda_{0}\frac{\mathbb{E}[u_{*}\boldsymbol{\varphi}_{+}'(c)]}{\mathbb{E}[u_{0}\boldsymbol{\varphi}_{+}'(c)]}+\alpha_{0},
\]
we obtain that $\mathbb{E}[\boldsymbol{\varphi}(c+\lambda_{0}u_{*})]=\mathbb{E}[\boldsymbol{\varphi}(\widetilde{c}+\alpha u_{0})]<\infty$,
which is a contradiction. Consequently, $\mathbb{E}[\boldsymbol{\varphi}(c+\lambda w_{*})]=\infty$
for all $\lambda>1$, and $w_{*}$ belongs to $\partial\mathcal{B}_{c}^{\varphi}$
and satisfies $\mathbb{E}[\boldsymbol{\varphi}(c+w_{*})]<\infty$. 

Proceeding as above, we show that $\mathbb{E}[\boldsymbol{\varphi}(c+\lambda w^{*})]<\infty$
for all $0\leq\lambda<1$, and $\mathbb{E}[\boldsymbol{\varphi}(c+\lambda w^{*})]=\infty$
for all $\lambda\geq1$. This result implies that $w^{*}$ belongs
to $\partial\mathcal{B}_{c}^{\varphi}$ and is such that $\mathbb{E}[\boldsymbol{\varphi}(c+w^{*})]=\infty$. 
\end{proof}

For a function $u$ in $\partial\mathcal{B}_{c}^{\varphi}$, the behavior
of the normalizing function $\psi(\lambda u)$ as $\lambda\uparrow1$
depends on whether $\boldsymbol{\varphi}(c+u)$ is $\mu$-integrable
or not. This behavior is partially elucidated by the following proposition:

\begin{prop}
\label{prop:behavior_psi} Let $u$ be a function in the boundary
of $\mathcal{B}_{c}^{\varphi}$. For $\lambda\in[0,1)$, denote $\psi_{u}(\lambda):=\psi(\lambda u)$,
whose right derivative we indicate by $(\psi_{u})_{+}'(\lambda)$.
If \textup{$\mathbb{E}[\boldsymbol{\varphi}(c+u)]<\infty$ }then $\psi_{u}(\lambda)=\psi(\lambda u)$
converges to some $\alpha\in(0,\infty)$ as $\lambda\uparrow1$\textup{.
}On the other hand, if \textup{$\mathbb{E}[\boldsymbol{\varphi}(c+u)]=\infty$
then} $(\psi_{u})_{+}'(\lambda)$ tends to $\infty$ as $\lambda\uparrow1$\textup{. }\end{prop}
\begin{proof}
Observing that the normalizing function $\psi$ is convex with $\psi(0)=0$,
we conclude that $\psi_{u}(\lambda)=\psi(\lambda u)$ is non-decreasing
and continuous in $[0,1)$. Moreover, $(\psi_{u})_{+}'(\lambda)$
is non-decreasing in $[0,1)$. Fix any function $u$ in the boundary
of $\mathcal{B}_{c}^{\varphi}$ such that $\mathbb{E}[\boldsymbol{\varphi}(c+u)]<\infty$.
Assume that $\psi(\lambda u)$ tends to $\infty$ as $\lambda\uparrow1$.
In this case, it is clear that 
\[
\boldsymbol{\varphi}(c+\lambda u-\psi(\lambda u)u_{0})\leq\boldsymbol{\varphi}(c+u\mathbf{1}_{\{u>0\}}-\psi(\lambda u)u_{0})\rightarrow0,\qquad\text{as }\lambda\uparrow1.
\]
Since $\boldsymbol{\varphi}(c+\lambda u-\psi(\lambda u)u_{0})\leq\boldsymbol{\varphi}(c+u\mathbf{1}_{\{u>0\}})$,
we can use the Dominated Convergence Theorem to write
\[
\mathbb{E}[\boldsymbol{\varphi}(c+\lambda u-\psi(\lambda u)u_{0})]\rightarrow0,\qquad\text{as }\lambda\uparrow1,
\]
which is a contradiction to $\mathbb{E}[\boldsymbol{\varphi}(c+\lambda u-\psi(\lambda u)u_{0})]=1$.
Thus $\psi(\lambda u)$ is bounded in $[0,1)$, and $\psi(\lambda u)$
converges to some $\alpha\in(0,\infty)$ as $\lambda\uparrow1$. 

Now consider any function $u$ in the boundary of $\mathcal{B}_{c}^{\varphi}$
satisfying $\mathbb{E}[\boldsymbol{\varphi}(c+u)]=\infty$. Suppose
that $(\psi_{u})_{+}'(\lambda)$ converges to some $\beta\in(0,\infty)$
as $\lambda\uparrow1$. Then $\psi_{u}(\lambda)=\psi(\lambda u)$
converges to some $\alpha\in(0,\infty)$ as $\lambda\uparrow1$. From
Fatou's Lemma, it follows that
\[
\mathbb{E}[\boldsymbol{\varphi}(c+u-\alpha u_{0})]\leq\liminf_{\lambda\uparrow1}\mathbb{E}[\boldsymbol{\varphi}(c+\lambda u-\psi(\lambda u)u_{0})]=1.
\]
Since $\varphi(t,\cdot)$ is convex, for any $\lambda\in(0,1)$, we
can write
\begin{align*}
\boldsymbol{\varphi}(c+\lambda u-\psi(\lambda u)u_{0}) & =\boldsymbol{\varphi}\Bigl(\lambda(c+u-\alpha u_{0})+(1-\lambda)\Bigl(c-\alpha u_{0}+\frac{\alpha-\psi(\lambda u)}{1-\lambda}u_{0}\Bigr)\Bigr)\\
 & \leq\lambda\boldsymbol{\varphi}(c+u-\alpha u_{0})+(1-\lambda)\boldsymbol{\varphi}\Bigl(c-\alpha u_{0}+\frac{\alpha-\psi(\lambda u)}{1-\lambda}u_{0}\Bigr).
\end{align*}
Observing that $\beta=\lim_{\lambda\uparrow1}(\psi_{u})_{+}'(\lambda)=\lim_{\lambda\uparrow1}[\alpha-\psi(\lambda u)]/(1-\lambda)$,
we can infer that
\[
\boldsymbol{\varphi}(c+\lambda u-\psi(\lambda u)u_{0})\leq\boldsymbol{\varphi}(c+u-\alpha u_{0})+\boldsymbol{\varphi}(c-\alpha u_{0}+\beta u_{0}),
\]
showing that $\boldsymbol{\varphi}(c+\lambda u-\psi(\lambda u)u_{0})$
is dominated by an integrable function. Thus, by the Dominated Convergence
Theorem, it follows that 
\[
\mathbb{E}[\boldsymbol{\varphi}(c+u-\alpha u_{0})]=\mathbb{E}[\lim_{\lambda\uparrow1}\boldsymbol{\varphi}(c+\lambda u-\psi(\lambda u)u_{0})]=\lim_{\lambda\uparrow1}\mathbb{E}[\boldsymbol{\varphi}(c+\lambda u-\psi(\lambda u)u_{0})]=1.
\]
The definition of $u_{0}$ tells us that $\mathbb{E}[\boldsymbol{\varphi}(\widetilde{c}+\lambda u_{0})]<\infty$
for all $\lambda\in\mathbb{R}$ and any measurable function $\widetilde{c}$
such that $\mathbb{E}[\boldsymbol{\varphi}(\widetilde{c})]=1$. In
particular, considering $\widetilde{c}=c+u-\alpha u_{0}$ and $\lambda=\alpha$,
we have that $\mathbb{E}[\boldsymbol{\varphi}(c+u)]<\infty$. This
contradicts the assumption that $\mathbb{E}[\boldsymbol{\varphi}(c+u)]=\infty$.
Therefore, $\lim_{\lambda\uparrow1}(\psi_{u})_{+}'(\lambda)=\infty$.
\end{proof}

\bibliographystyle{plain}
\bibliography{refs}

\end{document}